# Two presumptions in Gödel's interpretation of his own, formal, reasoning that are classically objectionable

## &

## Some significant consequences of interpreting arithmetical 'truth' constructively in Gödel's formal reasoning

Bhupinder Singh Anand[1]

Standard expositions of Gödel's 1931 paper on undecidable arithmetical propositions are based on two presumptions in Gödel's 1931 interpretation of his own, formal, reasoning - one each in Theorem VI and in Theorem XI - which do not meet Gödel's requirement of classically constructive, and intuitionistically unobjectionable, reasoning. We see how these objections can be addressed, and note some consequences thereof.

# I  Gödel's first, classically objectionable, presumption

## I-1  Introduction

The difficulty in situating an interpretation[2] of Gödel's first incompleteness theorem - Theorem VI in [Go31] - within the perspectives of classical (*pre-Cantorian*) reasoning,

---

[1] The author is an independent scholar. E-mail: re@alixcomsi.com; anandb@vsnl.com. Postal address: 32, Agarwal House, D Road, Churchgate, Mumbai - 400 020, INDIA. Tel: +91 (22) 2281 3353. Fax: +91 (22) 2209 5091.

[2] The word "interpretation" may be used both in its familiar, linguistic, sense, and in a mathematically precise sense; the appropriate meaning is usually obvious from the context.

arise simply because standard expositions[3] of Gödel's formal reasoning in his 1931 paper [Go31] on undecidable arithmetical propositions - which are, essentially, still based on Gödel's interpretation of his own, formal, reasoning in this paper - have yet to explicitly recognise that:

(*a*)   the truth of Gödel's unprovable proposition - which, syntactically, is of the form $[(\forall x)R(x)]$ - follows strictly from the axioms and rules of inference of Peano Arithmetic;

(*b*)   all the axioms and theorems of Peano Arithmetic are algorithmically (*Turing*) computable arithmetical truths (*when treated as Boolean functions*) under the standard interpretation;

(*c*)   there may be arithmetical relations - such as the one constructed by Gödel in his 1931 paper - which are meta-mathematically provable as instantiationally computable (*Turing-decidable*) truths, but which may not necessarily be algorithmically computable (*Turing-computable*) truths.

---

Mathematically, following Tarski ([Me64], §2, p49): "An *interpretation* consists of a non-empty set D, called the *domain* of the interpretation, and an assignment to each predicate letter $A_j{}^n$ of an *n*-place relation in D, to each function letter $f_j{}^n$ of an *n*-place operation in D (*i.e., a function from $D^n$ into D*), and to each individual constant $a_i$ of some fixed element of D. Given such an interpretation, variables are thought of as ranging over the set D, and ~, =>, and quantifiers are given their usual meaning. (*Remember that an n-place relation in D can be thought of as a subset of $D^n$, the set of all n-tuples of elements of D.*)"

We note that the interpreted relation $R'(x)$ is obtained from the formula $[R(x)]$ of a formal system P by replacing every primitive, undefined symbol of P in the formula $[R(x)]$ by an interpreted mathematical symbol. So the P-formula $[(\forall x)R(x)]$ interprets as the sentence $(\forall x)R'(x)$, and the P-formula $[\neg(\forall x)R(x)]$ as the sentence $\neg(\forall x)R'(x)$.

[3] Such as, for instance, [Be59], [Bo03], [Me64], [Ro87], [Sh67], [Sm92], [Wa64].

(*d*)    under a constructive interpretation of arithmetical 'truth' in Gödel's formal reasoning, any first-order Peano Arithmetic is naturally ω-inconsistent, and has no non-standard models (*one significant consequence of this is that P ≠ NP*[4]).

## I-2    The truth of Gödel's unprovable proposition follows strictly from the axioms and rules of inference of Peano Arithmetic

It is a common misconception that an arithmetical statement - such as the one constructed

by Gödel - can be true and yet not proven from the axioms and rules of inference of a

first order Peano Arithmetic.[5]

However, the truth of Gödel's unprovable proposition *does* follow strictly from the axioms and rules of inference of PA.

To see this, note that the first part of Gödel's argument in Theorem VI of his 1931 paper

[Go31] is, simply, that, if PA is consistent, then we can *mechanically* construct a PA

formula[6] - which, syntactically, is of the form $[(\forall x)R(x)]$ - such that:

(*i*)    the formula $[(\forall x)R(x)]$ - when expressed purely as a string of symbols - *does not*

*follow mechanically* from the axioms of PA as the last of any finite sequence of

PA-formulas, each of which is either a PA-axiom, or a consequence of one or

more of the formulas preceding it in the sequence, by the *mechanical*

application of the rules of inference of PA;

---

[4] For the significance of the PvsNP problem, see [Cook].

[5] This is the essence of the Gödelian argument initially put forward by Lucas [Lu61][Lu96], and defended vigorously by Penrose [Pe90][Pe94][Pe96].

[6] This follows from Gödel's Theorem VII in [Go31], and the fact that Gödel's formal system P is defined recursively, in the sense that there exists a Turing machine which will, when presented with any input string, halt and accept if the string is a well-formed formula of the system, and halt and reject otherwise.

(*ii*)  for any given numeral [*n*] - which 'represents' the natural number *n* in PA - each of the formulas [*R*(*n*)] - again, when expressed purely as a string of symbols - *does follow mechanically* from the axioms of PA as the last of some finite sequence of PA-formulas, each of which is either a PA-axiom, or a consequence of one or more of the formulas preceding it in the sequence, by the *mechanical* application of the rules of inference of PA.

Now, (*i*) is, essentially, the standard definition (*due to Gödel*) of the meta-assertion:

(*iii*)  The PA-formula [$(\forall x)R(x)$] is unprovable in PA,

whilst (*ii*) is, essentially, the standard definition (*due to Tarski*) of the meta-assertion:

(*iv*)  The PA-formula [$(\forall x)R(x)$] is true in PA.

Hence, by definition, the appropriate interpretation of Gödel's reasoning (*i*) and (*ii*) ought to be[7]:

(*v*)  The PA-formula [$(\forall x)R(x)$] is unprovable but true in PA.

However standard expositions of Gödel's formal reasoning in [Go31] assert only that:

(*vi*)  The PA-formula [$(\forall x)R(x)$] is unprovable in PA, but true in the standard interpretation of PA.

They, thus, fail to highlight the fact that (*i*) and (*ii*), both, actually, follow strictly from the axioms and rules of inference of PA, and that the meta-assertion:

(*vii*) The PA-formula [$(\forall x)R(x)$] is true in the standard interpretation of PA,

---

[7] This point is discussed in more detail in "Why we shouldn't fault Lucas and Penrose for continuing to believe in the Gödelian argument against computationalism": http://alixcomsi.com/Why_we_shouldnt.htm

classically, implies[8] the meta-assertion:

> (*viii*) The PA-formula [$R(x)$] is provable in PA whenever we substitute a numeral [$n$]
> for the variable [$x$] in [$R(x)$].

In other words, classically, (*vii*) also asserts the instantiational PA-provability of a denumerable infinity of arithmetical propositions.

The question arises: When is a PA-true formula, such as, say, [$(\forall x)R(x)$], also PA-provable?

## I-3    The axioms and theorems of Peano Arithmetic are algorithmically (*Turing*) computable truths under the standard interpretation

Now, it follows from Turing's 1936 paper [Tu36] that if, say, the formula [$R'(x_1, x_2, x_3 \ldots, x_n)$][9] is either an axiom, or a theorem[10], of a Peano Arithmetic (*such as the system P considered by Gödel in his 1931 paper*), then there is always a Turing-machine[11], T, such

---

[8] As we show later, the two are not equivalent. This is because the classical meta-statement, "The PA-formula [$(\forall x)R(x)$] is true in the standard interpretation of PA", is ambiguous, and can mean either of,

"The PA-formula [$(\forall x)R(x)$] is algorithmically decidable as true in the standard interpretation of PA", or,

"The PA-formula [$(\forall x)R(x)$] is instantiationally decidable as true in the standard interpretation of PA".

Note also, that the meta-statement, "The PA-formula [$\neg(\forall x)R(x)$] is true in PA", does not translate as the meta-statement, "It is not true that the PA-formula [$R(x)$] is provable in PA whenever we substitute a numeral [$n$] for the variable [$x$] in [$R(x)$]", but as the meta-statement, "It is not algorithmically decidable as true that the PA-formula [$R(x)$] is provable in PA whenever we substitute a numeral [$n$] for the variable [$x$] in [$R(x)$]".

[9] We use square brackets to indicate that the expression within the brackets is intended to be viewed purely as a syntactical string of mathematical symbols, devoid of any meaning whatsoever. In other words, it is not to be seen as representing any concept under an intuitive interpretation of the symbols.

[10] Formally, a theorem of a Peano Arithmetic, such as Gödel's system P, is a sequence of well-formed formulas of P, each of which is either an axiom, or an immediate consequence - by means of the rules of deduction of P - of any one, or more, of the preceding formulas of the sequence.

[11] As defined by Turing in his 1936 paper on computable numbers [Tu36]. For a more formal definition of algorithmic / Turing computability, see [Me64].

that, for any set of natural numbers $(a_1, a_2, a_3 \ldots, a_n)$, T will compute the arithmetical proposition $R'(a_1, a_2, a_3 \ldots, a_n)$ as TRUE in a finite number of steps (*Appendix A*).

Gödel's remarkable achievement in [Go31] lay in his actual construction - using only classically acceptable reasoning - of a formula of one variable in P, say $[R(x)]$, such that:

(*iii*)  $[(\forall x)R(x)]$ is not a theorem of P (*assuming that P is consistent*);

(*iv*)  for any natural number $a$, there is a Turing-machine $T_a$ that will compute the arithmetical proposition $R(a)$ as TRUE in a finite number of steps.

However, we cannot assume from this that there must be a single Turing machine T' such that, for any natural number $a$, T' will compute the arithmetical proposition $R(a)$ as TRUE in a finite number of steps.[12]

## I-4    The distinction between classical verifiability and Turing-computability

In other words, the precise definition of the algorithmic computability of number-theoretic functions - formalised by Turing in his 1936 paper [Tu36] on computable numbers - allows the following distinction to be made between the instantiational decidability of number-theoretic relations, and their algorithmic computability (*when treated as Boolean functions*):

**Definition 1**: A total number-theoretical relation, say $R'(x_1, x_2, ..., x_n)$, when treated as a Boolean function, is *Turing-decidable* in M if, and only if, it is instantiationally equivalent to a number-theoretic relation, $S(x_1, x_2, ..., x_n)$, and there is a Turing-machine T such that, for any given natural number sequence, $(a_1, a_2, ..., a_n)$, T will compute $S(a_1, a_2, ..., a_n)$ as either TRUE, or as FALSE.

---

[12] When interpreting Gödel's 1931 paper on undecidable arithmetical propositions [Go31], it is useful to keep in mind that this paper pre-dated, both, Church's development [Ch36] of the lambda calculus (*formalising and extending the primitive recursive arithmetic introduced by Gödel in this paper*), and Turing's 1936 paper [Tu36] on computable numbers, where Turing, essentially, introduced the concept of algorithmic computability.

**Definition 2**: A total number-theoretical relation, $R'(x_1, x_2, ..., x_n)$, when treated as a Boolean function, is *Turing-computable* in M if, and only if, there is a Turing-machine T such that, for any given natural number sequence, $(a_1, a_2, ..., a_n)$, T will compute $R'(a_1, a_2, ..., a_n)$ as either TRUE, or as FALSE.

Note that although Turing-decidability (*which is, essentially, verifiability in the classical sense*) of a number-theoretic relation necessarily implies Turing-computability, the converse need not hold.

Thus, classically, a number-theoretic relation, such as, say Gödel's $R(x)$, would be decidable (*Turing-decidable*) as true if it were instantiationally equivalent, for any given natural number, to a number-theoretic relation, $S(x)$, which could be interpreted as a Turing-algorithm, and which is Turing-computable as TRUE for all natural numbers.

However, we cannot conclude from this that $R(x)$ must also be interpretable as a Turing-algorithm, and be Turing-computable as TRUE for all natural numbers.

In other words - as in the case of Turing's Halting function - the definition of $R(x)$ may involve a, possibly implicit, circular reference to the totality of values of $R(x)$.

Consequently, any Turing-machine that computes $R(x)$ could go into a non-terminating loop for some input (*i.e., in the language of Turing-machines, an instantaneous tape description could be repeated in the course of the computation for that input*).

Now, in his 1931 paper [Go31], Gödel only shows that, if PA is consistent, then his $R(x)$ is Turing-decidable as always TRUE (*when treated as a Boolean function*).

The Turing-decidability of Gödel's $R(x)$ follows from the fact that $R(x)$ is defined constructively by Gödel (*as a consequence of Theorem VII of [Go31]*) so that it is instantiationally equivalent to a primitive recursive relation which, of course, is known to be Turing-computable.

However, the question remains:

Is $R(x)$ also Turing-computable as always TRUE (*when treated as a Boolean function*)?

If we interpret Gödel's reasoning with the above distinction in mind, then his PA-formula $[(\forall x)R(x)]$ may be[13] unprovable, but true, simply because the arithmetical relation $R(x)$ is a 'Halting' type of relation that is Turing-decidable as always TRUE, but not Turing-computable as always TRUE[14].

## I-5    Gödel's first, classically objectionable, presumption

It is intriguing to speculate on why the above distinction was not considered by Gödel even after Turing's publication of his paper [Tu36] on computable numbers.

---

[13] We show later that this is, indeed, so, since PA has no non-standard models, whence it follows that every total arithmetical relation that is Turing-computable as always TRUE is PA-provable.

[14] The reason, in this case, that $R(x)$ is Turing-decidable as always TRUE, but not Turing-computable as always TRUE, lies in the fact that Gödel's definition of $R(x)$ is in terms of a primitive recursive relation $Q(x)$, and involves a circularity in that the definition of $Q(x)$, implicitly, refers to the totality of values of $R(x)$.

As a consequence, there is some value of $n$ such that, when computing $R(n)$, any Turing-machine that computes $R(x)$ does not halt on the input $n$ since, as in the case of Turing's Halting function, an instantaneous tape description repeats itself.

A contributing factor may have been his 1931 commitment to the first of his classically objectionable presumptions: the assumption of ω-consistency for P will not lead to an inconsistency[15].

This allowed him to meta-mathematically argue further - again, using only classically acceptable reasoning - that the formula $[\neg(\forall x)R(x)]$ is, also, not a theorem of an ω-consistent P.

Another factor may have been Gödel's reluctant acceptance, along with Church, of the thesis (*now known as the Church-Turing Thesis*) identifying effectively computable[16] number-theoretic functions with algorithmically computable number-theoretic functions.

From the first, Gödel concluded his first incompleteness theorem (*Theorem VI of his 1931 paper*), to the effect that, if P can be assumed ω-consistent, then the formulas $[(\forall x)R(x)]$ and $[\neg(\forall x)R(x)]$ are, both, not provable in P.

---

[15] Gödel noted in his 1931 paper that, if P is ω-consistent, then it is necessarily consistent in the classical sense, although the converse need not be necessarily true. However, we show that this implication is false.

[16] Note that the concept of 'effective computability' is essentially intuitive.

It is noteworthy that Gödel (*initially*), and Church (*subsequently - not least because of Gödel's disquietitude*), enunciated Church's formulation of 'effective computability' as a Thesis, because they were instinctively uncomfortable with accepting it as a definition that fully captures the essence of '(*intuitive*) effective computability'.

Gödel's reservations seem vindicated if we accept that a number-theoretic function can be effectively computable instantiationally, but not algorithmically.

Since every algorithmically function is, necessarily, computable instantiationally, we can, then, indeed define a function as effectively computable (*intuitively*) if, and only if, it is computable instantiationally.

The essence of Church's Thesis would then be the postulation that every effectively computable function is (*instantiationally*) equivalent to a recursive function.

Note that, in its standard form, Church's Thesis postulates a strong identity - and not simply a weak equivalence - between an effectively computable function and some recursive function.

The second may, quite possibly, have prevented him from recognising the possibility, and significance, of number-theoretic relations that are effectively computable as true instantiationally, but not algorithmically.

## I-6   Why the assumption of ω-consistency in a Peano Arithmetic is classically objectionable, and may invite inconsistency

It is interesting to see how, and why, Gödel's presumption - that P can be assumed ω-consistent without inviting inconsistency - obscures the significance of his primary achievement - that of constructing a number-theoretic, arithmetical, relation that may be effectively computable as true instantiationally, but not necessarily algorithmically.

Now, Gödel defines P as ω-consistent if, and only if, there is no P-formula, such as, say, his [$R(x)$], such that [¬($\forall x)R(x)$] is P-provable, and also that [$R(n)$] is P-provable for any numeral [$n$] of P.

The possible inconsistency in an ω-consistent P surfaces only when we try to interpret the above definition verifiably in our classical arithmetic based on Dedekind's Peano Postulates (*termed as the "standard interpretation of Peano Arithmetic"*), which Gödel's system P was intended to formalise recursively (*hence, implicitly, algorithmically*).

Now, it is reasonable to suppose that Gödel's tacit justification for the presumption of ω-consistency as an intuitively natural, and consistent, meta-property of P lay in the following argument:

(*a*)  if [¬($\forall x)R(x)$] is P-provable, then the meta-proposition, "*R(x)* is true for all natural numbers $x$", is false,

and:

(*b*) if [*R*(*n*)] is P-provable for any numeral [*n*] of P, then the meta-proposition, "*R*(*x*) is true for all natural numbers *x*", is true.

However, if we note, first, that the phrase, "*R*(*x*) is true for all natural numbers *x*", is ambiguous, and can be taken to mean either, "*R*(*x*) is effectively computable as true algorithmically (*Turing-computable*) for all natural numbers *x*", or, "*R*(*x*) is effectively computable as true instantiationally (*Turing-decidable*) for all natural numbers *x*", and, second, that the latter need not, necessarily, imply the former, then we can, indeed, have that:

(*c*) if [¬(∀*x*)*R*(*x*)] is P-provable, then the meta-proposition, "*R*(*x*) is effectively computable as true algorithmically (*Turing-computable*) for all natural numbers *x*", is false

and, at the same time, that:

(*d*) if [*R*(*n*)] is P-provable for any numeral [*n*] of P, then the meta-proposition, "*R*(*x*) is effectively computable as true instantiationally (*Turing-decidable*) for all natural numbers *x*", is true.

## I-7 A constructive interpretation of ω-consistency implies that Gödel's (*first incompleteness*) Theorem VI may be vacuously true

It follows that Gödel's presumption that assuming ω-consistency for P will not lead to an inconsistency is questionable, and hence his argument for the construction of an undecidable proposition in P, cannot be considered classically constructive, and intuitionistically unobjectionable.

It also follows that we cannot, further, conclude from Gödel's reasoning that his proposition, [¬(∀*x*)*R*(*x*)], is necessarily not P-provable.

## I-8  Rosser's argument, too, implicitly assumes that Peano Arithmetic is ω-consistent

Contrary to the accepted view, Rosser's 'extension' of Gödel's argument [Ro36] also assumes - albeit, implicitly - that P is ω-consistent.

Thus, in his proof, Rosser argues at one point[17] that if we assume that a formula such as, say, $[(\exists x)F(x)]$, is provable in a consistent P, then we can assume, further, that there is some P-formula, say $[a]$, such that $[F(a)]$ is provable.

He, then, argues that the P-provability of $[F(a)]$ leads to a contradiction, and concludes that $[(\exists x)F(x)]$ cannot be P-provable.

Now, if P were ω-inconsistent, then, since, formally, $[(\exists x)F(x)]$ is defined as $[\neg(\forall x)\neg F(x)]$, we can, indeed, have, both, that $[\neg(\forall x)\neg F(x)]$ is P-provable, and that $[\neg F(x)]$ is also P-provable for all substitutions of a P-formula $[a]$ for $[x]$.

Hence, unless we assume that P is ω-consistent, we cannot conclude - as Rosser does - that, from the P-provability of $[\neg(\forall x)\neg F(x)]$, there is, necessarily, some P-formula, say $[a]$, such that $[F(a)]$ is provable.[18]

## I-9  If Peano Arithmetic is ω-inconsistent, then it may have no non-constructive, non-standard, models

Now, if we cannot conclude that $[\neg(\forall x)R(x)]$ is necessarily not P-provable, then we cannot conclude that we can add $[\neg(\forall x)R(x)]$ as an axiom to P and obtain a consistent,

---

[17] I analyse a standard exposition of Rosser's proof in: Reviewing Gödel's and Rosser's meta-reasoning of undecidability (*http://alixcomsi.com/Constructivity_consider.htm*).

[18] I argue that this, essentially, lay at the heart of Brouwer's objection to Hilbert's unrestricted use of the Law of the Excluded Middle. See: Why Brouwer was right in suggesting that Hilbert's Law of the Excluded Middle needed qualification (*http://alixcomsi.com/Why_Brouwer_was_right.htm*).

extended, system of Peano Arithmetic, which contains numbers that are not natural numbers.

However, it is precisely the unvalidated presumption by Gödel that P+[¬(∀x)R(x)] is consistent [Go31], which appears to underpin the postulation of non-standard models[19] of Peano Arithmetic.

These are models that postulate the existence of numbers, satisfying the axioms of Peano Arithmetic, which exist in the domain of an interpretation, and can, possibly, be named, but which are not numerals, and, so, cannot be constructed within the Arithmetic.

In Appendix B we show, however, that all the elements in the domain of any model of PA are necessarily the numerals (*which, by definition, are the successors of 0*) and so PA has no non-standard models.

So, if we cannot add [¬(∀x)R(x)] as an axiom to P without inviting inconsistency, it follows that P is ω-inconsistent.

## I-10   A significant consequence of a constructive interpretation of Gödel's reasoning in Theorem VI of his 1931 paper: P ≠ NP

We are, thus, compelled to seek alternative, and more constructive, interpretations of the significance, and consequences, of Gödel's reasoning in his 1931 paper.[20]

---

[19] We define [cf. [Me64], p49-53] a model as an interpretation of a set of well-formed formulas if, and only if, every well-formed formula of the set is true for the interpretation (*by Tarski's definitions of the satisfiability and truth of the formulas of a formal language under an interpretation*).

We define a non-standard model, say M', of first-order Peano Arithmetic as one in which, if [(∀x)R(x)] is Gödel's undecidable formula, then ~R(s) holds for some $s$ in the domain of M' that is not a natural number (*i.e., $s$ is not a successor of 0*).

[20] See also: PA is instantiationally and arithmetically complete, but algorithmically incomplete: An alternative interpretation of Gödelian incompleteness under Church's Thesis that links formal logic and computability (*http://alixcomsi.com/PA_is_instantiationally_complete.htm*).

One of the more significant consequences is that if a total arithmetical relation is Turing-computable as always TRUE, then it is PA-provable.

For, if we assume that there is a total arithmetical relation that is Turing-computable as always TRUE, but which is not the standard interpretation of a PA-provable formula, then this implies that there is a non-standard model of PA.

Another follows from the argument that (cf. [Ra02]), if there is no polynomial time algorithm $A$ that gets as input a Boolean formula $f$ and outputs 1 if and only if $f$ is a tautology, then $P \neq NP$.

Hence Gödel's reasoning in Theorem VI of [Go31] - to the effect that there is a total arithmetical relation that is Turing-decidable as always TRUE, but which is not Turing-computable as always TRUE - implies that $P \neq NP$.

## II    Gödel's second, classically objectionable, presumption

### II-1    Gödel's reasoning in his (second incompleteness) Theorem XI

The second of Gödel's classically objectionable presumptions pertains to the reasoning by which he constructs another P-formula, [*W*] (*say, with Gödel-number w*), and concludes that, when [*W*] is interpreted in our intuitive arithmetic of the natural numbers (*as covered by Dedekind's formulations of the Peano Axioms*), it intuitively translates as an arithmetical proposition that is true if, and only if, a specified formula of P is unprovable in P.

Now, if there were such a P-formula as [*W*] which interprets, and intuitively translates, as above, and we could prove [*W*] in P, then the above would imply that there is a proof sequence in P to the effect that P contains an unprovable formula, and so we would have proved that P is consistent by a proof sequence within P.

This follows since, classically, an inconsistent system necessarily proves every well-formed formula of the system.

Now Gödel shows - again by classical reasoning - that his formula [*W*] is not P-provable.

From this he concludes his second incompleteness theorem (*Theorem XI of his 1931 paper*), to the effect that the consistency of a Peano Arithmetic cannot be proven within the Arithmetic.

## II-2    Gödel presumes that he has constructed an arithmetical proposition that translates intuitively as true if, and only if, P is consistent

Now, the classically objectionable presumption made by Gödel is that, when the P-formula [*W*] is interpreted in our intuitive arithmetic of the natural numbers (*as covered by Dedekind's formulations of the Peano Axioms*), it intuitively translates as an arithmetical proposition, say *W*, that is true if, and only if, a specified formula of P is unprovable in P.

To see the classically objectionable step of Gödel's reasoning in Theorem XI of his 1931 paper, it is necessary to see how Gödel shows, in a constructive, and intuitionistically unobjectionable manner, how meta-propositions of a formal Peano Arithmetic, say P, can be defined by means of primitive recursive functions and relations using only classical reasoning.

## II-3    Gödel defines 45 number-theoretic primitive recursive functions and relations

Amongst the 46 such number-theoretic functions and relations defined by Gödel, in a classically constructive manner, are the following:

(#12) A number-theoretic primitive recursive relation, *Var(x)*, which is true if, and only if, $x$ is the Gödel-number[21] of a variable of P.

(#13) A number-theoretic primitive recursive function, *Neg(x)*, which is the Gödel-number of the negation of the formula of P whose Gödel-number is $x$.

(#14) A number-theoretic primitive recursive relation, *xDisy*, which is true if, and only if, the disjunction of the formulas of P, whose Gödel-numbers are $x$ and $y$, is true.

(#15) A number-theoretic primitive recursive function, *xGeny*, which is true if, and only if, the generalisation of the formula of P whose Gödel-number is $y$, by the the formula whose Gödel-number is $x$, is true (*assuming that x is the Gödel-number of a formula of P that is a variable*).

(#17) A number-theoretic primitive recursive function, *Z(n)*, which is the Gödel-number of the numeral of P that represents the natural number $n$.

(#20) A number-theoretic primitive recursive relation, *Elf(x)*, which is true if, and only if, $x$ is the Gödel-number of an elementary formula of P.

(#22) A number-theoretic primitive recursive relation, *FR(x)*, which is true if, and only if, $x$ is the Gödel-number of a sequence of formulas of P, each one of which is either an elementary formula of P, or comes from the preceding ones by the operations of negation, disjunction, or generalisation.

---

[21] In [Go3], p13, Gödel sets up a 1-1 correspondence of natural numbers to the primitive symbols, formulas, and finite sequences of formulas, of PA, in such a manner that various meta-mathematical functions and relations involving concepts such as, e.g., "variable", "formula", "sentence", "axiom", "provable formula", etc., can be expressed as equivalent to number-theoretic functions and relations.

(#23) A number-theoretic primitive recursive relation, *Form(x)*, which is true if, and only if, *x* is the Gödel-number of a formula of P (*i.e., the last term of a sequence of formulas*).

(#26) A number-theoretic primitive recursive relation, *vFrx*, which is true if, and only if, *x* is the Gödel-number of a formula of P in which the variable of P, whose Gödel-number is *v*, is free.

(#29) A number-theoretic primitive recursive function, A($v$, $x$), which is the number of places at which the variable of P, whose Gödel-number is *v*, is free in the formula of P whose Gödel-number is *x*.

(#42) A number-theoretic primitive recursive relation, *Ax(x)*, which is true if, and only if, *x* is the Gödel-number of an axiom of P.

(#43) A number-theoretic primitive recursive relation, *Fl($x$, $y$, $z$)*, which is true if, and only if, *x* is the Gödel-number of a formula of P that is an immediate consequence of the formulas of P whose Gödel-numbers are *y* and *z*.

(#44) A number-theoretic primitive recursive relation, *Bw(x)*, which is true if, and only if, *x* is the Gödel-number of a proof sequence in P (*a finite sequence of formulas of P each of which is either an axiom of P, or an immediate consequence of two preceding formulas of the sequence*).

(#45) A number-theoretic primitive recursive relation, *xBy*, which is true if, and only if, *x* is the Gödel-number of a proof sequence of P whose last formula has the Gödel-number *y*.

The classically constructive, and intuitionistically unobjectionable, nature of the above, 45, representation of various meta-mathematical propositions by primitive recursive functions and relations is assured by Gödel by the following specification (*in a footnote*):

> "Everywhere in the following definitions where one of the expressions '$(\forall x)$', '$(\exists x)$', '(There is a unique $x$)' occurs it is followed by a bound for $x$. This bound serves only to assure the recursive nature of the defined concept."

In sharp contrast to the above 45 definitions, Gödel explicitly defines one further meta-mathematical relation that cannot be asserted to be recursive:

> (#46) The number-theoretic relation, *Bew*($x$), which is true if, and only if, $x$ is the Gödel-number of a provable formula of P.

Gödel's actual definition of (#46) is:

> (#46) The number-theoretic primitive recursive relation, *Bew*($x$), is true if, and only if, the number-theoretic relation $(\exists y)y\mathrm{B}x$ is true.

The significance of the distinction, between the first 45 definitions and the 46[th], is that, by Gödel's Theorem VII (*which is independent of his Theorem VI, and should actually have preceded it*), it can be shown that any recursive relation, say $Q(x)$, can be represented in P by some, corresponding, arithmetical formula, say $[R(x)]$, such that, for any natural number $n$:

> If $Q(n)$ is true, then $[R(n)]$ is P-provable

> If $Q(n)$ is false, then $[\neg R(n)]$ is P-provable.

Now, Gödel's Theorem VI establishes that the above representation does not extend to the closure of a recursive relation, in the sense that we cannot assume:

If $(\forall x)Q(x)$ is true (*over all the natural numbers*), then $[(\forall x)R(x)]$ is P-provable.

In other words, we cannot assume that, even though the recursive relation, $Q(x)$, can be represented in P by some arithmetical formula, the number-theoretic proposition $(\forall x)Q(x)$ must, necessarily, also be representable similarly.

However this, precisely, is the presumption made by Gödel in his Theorem XI (*where he concludes that the consistency of P is not P-provable*).

Thus, he first defines the notion of "P is consistent" in a classically constructive, and intuitionistically unobjectionable, manner, as follows:

P is consistent if, and only if, the number-theoretic relation *Wid*(P) is true.

Here, *Wid*(P) is defined by the number-theoretic relation:

$(\exists x)[Form(x)$ & $\neg(Bew(x)]$

The above translates intuitively as:

There is a natural number $x$ such that $x$ is the Gödel-number of a formula of P, and this formula is not P-provable.

That this implies the consistency of P follows from the classical definition of the consistency of any formal axiomatic theory, since, if the theory were inconsistent, then every well-formed formula of the theory would be provable within the theory.

Now, Gödel's classically objectionable step lies in his presumption that *Wid*(P) can be represented (*[Go31] uses the term "expressed"*) by some formula [*W*] of P, with Gödel-number $w$, that, under the standard interpretation, maintains its intuitive meaning.

To see why the presumption is objectionable, it is sufficient to note that - again without adequate justification - Gödel also presumes in the proof of his Theorem XI [Go31] that if the recursive relation, $Q(x, p)$, is represented by the P-formula $[R(x, p)]$, whose Gödel-number is $r$, then $(\forall x)Q(x, p)$ can be "expressed" in P by the formula $[(\forall x)R(x, p)]$, whose Gödel-number is $17Genr$, in such a manner that the proposition "$[(\forall x)R(x, p)]$ is true" has the same intuitive meaning, under the standard interpretation, as "$(\forall x)Q(x, p)$ is true".

As noted earlier, such an assumption does, indeed, follow if $R(x, p)$ is Turing-computable as always TRUE, but it may not follow if $R(x, p)$ is Turing-decidable as always TRUE, but not Turing-computable as always TRUE.

In other words, if $[W]$, too, is unprovable but true in P, then the consistency of P is, in fact, provable instantiationally in P.

Hence, at best, Gödel's reasoning can only be taken to establish that the consistency of P is unprovable algorithmically in P.

The significance of this is that we can no longer argue that if the consistency of a system that contains a Peano Arithmetic is, somehow, provable within the system, then the system must be inconsistent.

This conclusion can only be drawn for such a system if it purports to prove its own consistency algorithmically.

## Appendix A: The Theorems of PA interpret as Turing-computable truths

Consider the axioms of PA:

(S1)  $[(x_1 = x_2) \rightarrow ((x_1 = x_3) \rightarrow (x_2 = x_3))]$;

(S2)  $[(x_1 = x_2) \rightarrow (x_1' = x_2')]$;

(S3)  $[0 \neq (x_1)']$;

(S4)  $[((x_1)' = (x_2)') \rightarrow (x_1 = x_2)]$;

(S5)  $[(x_1 + 0) = x_1]$;

(S6)  $[(x_1 + x_2') = (x_1 + x_2)']$;

(S7)  $[(x_1 * 0) = 0]$

(S8)  $[(x_1 * (x_2')) = ((x_1 * x_2) + x_1)]$;

(S9)  For any well-formed formula $[F(x)]$ of PA:

$[F(0) \rightarrow ((\forall x)(F(x) \rightarrow F(x')) \rightarrow (\forall x)F(x))]$.

Under the standard interpretation, M, each of the PA axioms - except Induction - is Turing-computable as always TRUE in the following sense:

> **Definition**: A total number-theoretical relation, $R(x_1, x_2, ..., x_n)$, when treated as a Boolean function, is Turing-computable in M if, and only if, there is a Turing-machine T such that, for any given natural number sequence, $(a_1, a_2, ..., a_n)$, T will compute $R(a_1, a_2, ..., a_n)$ as either TRUE, or as FALSE.

(*i*) If $[R]$ is an atomic formula $[R(a_1, a_2, ..., a_n)]$ of PA, then $[R]$ is Turing-computable as TRUE for the natural number input $(a_1, a_2, ..., a_n)$ if, and only if, the arithmetical relation $R(a_1, a_2, ..., a_n)$ is Turing-computable as TRUE on the natural number input $(a_1, a_2, ..., a_n)$;

(*ii*) The PA-formula [¬$R$] is Turing-computable as TRUE for the natural number input ($a_1, a_2, ..., a_n$) if, and only if, [$R$] is Turing-computable as FALSE for the natural number input ($a_1, a_2, ..., a_n$);

(*iii*) The PA-formula [$R => S$] is Turing-computable as TRUE for the natural number input ($a_1, a_2, ..., a_n$) if, and only if, either [$R$] is Turing-computable as FALSE for the natural number input ($a_1, a_2, ..., a_n$), or [$S$] is Turing-computable as TRUE for the natural number input ($a_1, a_2, ..., a_n$);

(*iv*) The PA-formula [($\forall x_i$)$R$] is Turing-computable as TRUE for the natural number input ($a_1, a_2, ..., a_n$) if, and only if, [$R$] is Turing-computable as TRUE for every natural number input ($b_1, b_2, ..., b_n$) where ($b_1, b_2, ..., b_n$) differs from ($a_1, a_2, ..., a_n$) in at most the $i$'th component.

Further, it is easily seen that the following rules of Inference in PA preserve the above Turing-computability under the standard interpretation:

(*vii*) Modus Ponens: [$B$] follows from [$A$] and [$A \rightarrow B$];

(*viii*) Generalisation: [($\forall x$)$A$] follows from [$A$].

It follows that the Theorems of PA are Turing-computable as always TRUE in the above sense.

# Appendix B: Standard, first-order, PA has no non-standard model

We give an elementary proof that, if PA is a standard, first-order, Peano Arithmetic - as defined in Appendix A - then PA has no non-standard model (*i.e., a model whose domain contains an element that is not a successor of 0*).

We denote by *G*(*x*) the PA-formula:

    [*x*=0 v ¬(∀*y*)¬(*x*=*y*')].

This translates, under every interpretation of PA, as:

    Either *x* is 0, or *x* is a 'successor'.

Now, in every interpretation of PA, we have that:

    (*a*) *G*(0) is true;

    (*b*) If *G*(*x*) is true, then *G*(*x*') is true.

It follows, from Gödel's completeness theorem[22], that:

    (*a*) [*G*(0)] is provable in PA;

    (*b*) [*G*(*x*) => *G*(*x*')] is provable in PA.

We also have, by Generalisation (*Appendix A*), that:

    (*c*) [(∀*x*)(*G*(*x*) => *G*(*x*'))] is provable in PA;

From the Induction axiom S9 (*Appendix A*), we thus have that:

    (*d*) [(∀*x*)*G*(*x*)] is provable in PA.

We conclude that, except 0, every element in the domain of any interpretation of PA is a successor of 0.

Since, by definition, the successors of 0 are the natural numbers, it follows that there are no non-standard models of a standard, first-order, Peano Arithmetic.

---

[22] Cf. [Me64], p68.